\documentclass[12pt]{amsart}

\newtheorem{theorem}{Theorem}
\newtheorem{lemma}[theorem]{Lemma}

\newtheorem{corollary}[theorem]{Corollary}

\theoremstyle{definition}
\newtheorem{example}{Example}

\newcommand{\ltwo}{L^2({\mathbb R})}
\newcommand{\bbz}{\mathbb{Z}}
\newcommand{\bbrn}{\mathbb{R}^n}
\newcommand{\Cal}{\mathcal}
\newcommand{\Bb}{\mathbb}

\newcommand{\wo}{\widehat{W_0}}

\begin{document}

% MAKE THE TITLE
\title{On the Translation Invariance of Wavelet Subspaces}
\author{Eric Weber}
\address{Department of Mathematics, University of Colorado,
Boulder, CO 80309-0395} 
\curraddr{Department of Mathematics, Texas A\&M University, College Station, TX 77843-3368}
\email{weber@math.tamu.edu}
\begin{abstract}
An examination of the translation invariance of $V_0$ under dyadic rationals is presented, generating a new equivalence relation on the collection of wavelets.  The equivalence classes under this relation are completely characterized in terms of the support of the Fourier transform of the wavelet.  Using operator interpolation, it is shown that several equivalence classes are non-empty.
\end{abstract}
\maketitle

% REMEMBER TO DOUBLE SPACE BEFORE SENDING TO JFAA

% \setlength{\baselineskip}{30pt}

% BEGIN THE MAIN BODY OF THE PAPER
\section{Introduction}

A wavelet $\psi \in \ltwo$ is a complete wandering vector for the unitary system $\{D^n T^l: n, \ l \in \Bb{Z} \}$, i.e. the collection $\{ D^n T^l \psi : n, \ l \in \Bb{Z} \}$ is an orthonormal basis for $\ltwo$, where $D$, $T$ are defined on $\ltwo$ as: $Df(x) = \sqrt{2} f(2x)$ and $Tf(x) = f(x - 1)$.  Every wavelet can be associated with a \emph{Generalized Multiresolution Analysis}, or GMRA (see \cite{BCMO}).  Indeed, define the subspaces $V_j = \overline{span} \{D^n T^l \psi : n < j,\ l \in \bbz\}$, then it is routine to verify that these subspaces satisfy the following four conditions:
\begin{enumerate}
\item $V_j \subset V_{j+1}$,
\item $DV_j = V_{j+1}$,
\item $\cap_{j \in \bbz} V_j = \{0\}$ and $\cup_{j \in \bbz} V_j$ has dense span in $\ltwo$,
\item $V_0$ is invariant under $T$. \label{I:VIT}
\end{enumerate}
We shall call $V_0$ the \emph{core space} for $\psi$.  Item \ref{I:VIT} is of interest since the core space is invariant under translations by the group $\Bb{Z}$.  A natural question is:  are there other groups of translations under which $V_0$ is invariant?  This paper will answer this question by looking at groups of translations by dyadic rationals.

Denote by $T_{\alpha}$ the unitary operator $T_{\alpha}f(x) = f(x - \alpha)$.  $T$ is to be understood as $T_1$.  Note that $\widehat{T_{\alpha}} = M_{e^{-i \alpha \xi}}$.  In this paper, we shall consider the groups of translations $\Cal{G}_n = \{ T_{\frac{m}{2^n}}: m \in \Bb{Z} \}$, and the group $\Cal{G}_{\infty} = \{ T_{\alpha} : \alpha \in \Bb{R} \}$.  Denote by $\Cal{L}_n$ the collection of all wavelets whose core space is invariant under $\Cal{G}_n$.  Note that these collections are nested:
\[
\Cal{L}_0 \supset \Cal{L}_1 \supset \Cal{L}_2 \ldots \supset \Cal{L}_n \supset \Cal{L}_{n+1} \supset \ldots \supset \Cal{L}_{\infty}
\]

We can then define an equivalence relation whose equivalence classes are given by $\Cal{M}_n = \Cal{L}_n - \Cal{L}_{n+1}$, with $\Cal{M}_{\infty} = \Cal{L}_{\infty}$.  Hence, $\Cal{M}_n$ is the collection of all wavelets such that $V_0$ is invariant under $\Cal{G}_n$ but not $\Cal{G}_{n+1}$.  The goal of this paper is to characterize these equivalence classes, while showing that several of them are not empty.

In general, $V_0$ can be quite complicated in structure.  Indeed, it may not even be generated by translations of a finite number of functions. Hence, we wish to restrict our analysis to $W_0$.  Recall that $W_j$ is defined by $V_{j+1} = V_j \oplus W_j$.  Clearly, $W_j = \overline{span} \{D^j T^l \psi: l \in \Bb{Z} \}$.

\begin{lemma} \label{L:VG}
Let $r < n$ be integers, and let $p = n-r$.  The space $V_r$ (resp. $W_r$) is invariant under $\Cal{G}_k$ if and only if the space $V_n$ (resp. $W_n$) is invariant under $\Cal{G}_{k+p}$.
\end{lemma}
\begin{proof}
By definition, $f \in V_r$ if and only if $D^p f \in V_n$.  Suppose that $g \in V_n$ and define $f \in V_r$ such that $D^pf = g$.  Consider the following commutation relation:
\begin{align*}
D^p T_{\frac{m}{2^k}} f(x) &= 2^{\frac{p}{2}} T_{\frac{m}{2^k}} f(2^p x) \\
	&= 2^{\frac{p}{2}} f(2^p x + \frac{m}{2^k}) \\
	&= D^p f(x + \frac{m}{2^{k+p}}) \\
	&= T_{\frac{m}{2^{k+p}}} D^p f(x).
\end{align*}
This calculation establishes the statement.
\end{proof}

By lemma \ref{L:VG}, another way to describe $\Cal{M}_n$ is that $\psi \in \Cal{M}_n$ if $n$ is the largest integer such that $V_{-n}$ is invariant under integral translations.  If $\psi \in \Cal{M}_n$, we shall say $\psi$ has the \emph{translation invariance of order n} property.

\begin{theorem} \label{T:VW}
The core space $V_0$ for $\psi$ is invariant under the action of $\Cal{G}_n$ if and only if $W_0$ is invariant under the action of $\Cal{G}_n$.
\end{theorem}
\begin{proof}[If]  
Suppose that $W_0$ is invariant under $\Cal{G}_n$.  Then, by lemma~\ref{L:VG}, for $k > 0$, $W_k$ is invariant under $\Cal{G}_{n+k}$, whence $V_0^{\perp} = \oplus_{k=0}^{\infty} W_k$ is invariant under $\Cal{G}_n$.  If follows that $V_0$ is invariant under $\Cal{G}_n$.

\medskip
\noindent \emph{Only If.}
Suppose $V_0$ is invariant under $\Cal{G}_n$.  Then, again by lemma~\ref{L:VG}, $V_1$ is invariant under $\Cal{G}_{n+1}$, and hence $\Cal{G}_n$.  Since $V_1 = V_0 \oplus W_0$, it follows that $W_0$ is also invariant under $\Cal{G}_n$.
\end{proof}

For the purposes of this paper, we shall say that a set $E \subset \Bb{R}$ is \emph{partially self-similar} with respect to $\alpha \in \Bb{R}$ if there exists a set $F$ of non-zero measure such that both $F$ and $F + \alpha$ are subsets of $E$.  Additionally, if $G,H$ are two subsets of $\Bb{R}$, we shall say that $G$ is $2 \pi$ translation congruent to $H$ if there exists a measurable partition $G_n$ of $G$ such that the collection $\{G_n + 2 n \pi: n \in \Bb{Z}\}$ forms a partition of $H$, modulo sets of measure zero.  The letter $\lambda$ will denote Lebesgue measure.  Define a mapping $\tau: \Bb{R} \to [0, 2\pi)$ such that $\tau(x) - x = 2 \pi k$ for some integer $k$.

\section{A Characterization of $\Cal{M}_{\infty}$}

Recall that a wavelet set is a set $W \subset \mathbb{R}$ such that the function $\psi$ defined by $\hat\psi = \frac{1}{\sqrt{2 \pi}} \chi_W$ is a wavelet.  Such a wavelet $\psi$ is called a Minimally Supported Frequency (MSF) wavelet.  The following theorem reveals some structure of wavelet sets.
\begin{theorem} \label{T:ef}
Let $f \in \ltwo$, and let $E = supp(f)$.  Then $\{e^{-i n x} f(x) : n \in \Bb{Z}\}$ is an orthonormal basis for $L^2(E)$ if and only if the following two conditions hold:
\begin{enumerate}
\item supp(f) is $2\pi$ translation congruent to $[0, 2\pi)$,
\item $|f(x)| = \frac{1}{\sqrt{2\pi}}$ for almost every $x$.
\end{enumerate}
\end{theorem}

This next theorem completely characterizes wavelet sets; the proof of which uses the previous theorem.
\begin{theorem} \label{T:set}
Let $W \subset \Bb{R}$.  Then $W$ is a wavelet set if and only if the following two conditions hold:
\begin{enumerate}
\item $W$ is $2\pi$ translation congruent to $[0, 2\pi)$,
\item $\cup_{j \in \Bb{Z}} 2^j W = \Bb{R}$
\end{enumerate}
modulo null sets.
\end{theorem}
The proof of both theorems can be found in \cite{DL}, chapter 4.

\begin{theorem} \label{T:MSF}
Let $\psi$ be a wavelet.  Then, the following are equivalent:
\begin{enumerate}
\item[i)] $\psi$ is a MSF wavelet,
\item[ii)] the subspace $V_0$ is invariant under translations by all real numbers,
\item[iii)] the subspaces $V_j$ of the corresponding GMRA are invariant under integral translations.
\item[iv)] the subspaces $W_j$ of the corresponding GMRA are invariant under integral translations.
\end{enumerate}
\end{theorem}

\begin{proof} i) $\Rightarrow$ ii).
If $\psi$ is a MSF wavelet with wavelet set $W$, then by theorem~\ref{T:ef}, $\widehat{W_0} = L^2(W)$.  Clearly, $\forall \alpha \in \Bb{R}$, $W_0$ is invariant under $T_{\alpha}$ since $\widehat{W_0}$ is invariant under multiplication by $e^{-i \alpha \cdot}$.  It follows by theorem~\ref{T:VW} that $V_0$ is invariant under all translations.

\medskip
\noindent ii) $\Rightarrow$ iii).
Since $V_0$ is invariant under $\Cal{G}_n$ for all $n \geq 0$, by lemma~\ref{L:VG}, $V_{-n}$ is invariant under $\Cal{G}_0$.

\medskip
\noindent iii) $\Rightarrow$ iv).
By definition, $V_{j+1} = V_j \oplus W_j$.  If both $V_{j+1}$ and $V_j$ are invariant under integral translations, it follows immediately that $W_j$ is also invariant under integral translations.

\medskip
\noindent iv) $\Rightarrow$ i).
Let $\Cal{C}$ be the collection of all operators for which $W_0$ is invariant.  An easy calculation shows that $\Cal{C}$ is WOT (weak operator topology) closed.

If $W_j$ is invariant under integral translations, then again by lemma~\ref{L:VG}, $W_0$ is invariant under $\Cal{G}_j$ for all $j$.  Since $\cup_{n \geq 0} \Cal{G}_n$ is dense in $\Cal{G}_{\infty}$, in the strong operator topology, it follows that $W_0$ is invariant under $\Cal{G}_{\infty}$.  If we take the Fourier transform, then we get that $\wo$ is invariant under multiplication by $e^{-i \alpha \xi}$.  The linear span of these operators are dense in the collection $\{ M_h : h \in L^{\infty}(\Bb{R}) \}$ with respect to the WOT.  It follows that $\wo$ is invariant under multiplication by any $L^{\infty}(\Bb{R})$ function.

Next, we wish to show that $\wo = L^2(E)$, where $E = supp(\hat{\psi})$.  First note that since $\{ e^{-in \xi} \hat{\psi}(\xi)\}$ forms an orthonormal basis for $\wo$, $\hat{\psi}(\xi)$ has maximal support in the sense that if $\hat{f} \in \wo$, then the support of $\hat{f}$ is contained in the support of $\hat{\psi}$.  This immediately implies that $\wo \subset L^2(E)$.

Let $g(\xi)$ be a compactly supported simple function, whose support $F$ is contained in $E$.  Define $E_n = \{\xi : \frac{1}{n-1} \geq \hat{\psi}(\xi) > \frac{1}{n}\}$, and define $F_n = F \cap E_n$.  Since $g$ is a simple function, it is uniformly bounded by some constant $M$.  Let $\epsilon > 0$ be given.  Choose an $N$ such that $\lambda(\cup_{n>N} F_n) < \frac{\epsilon}{M}$, and define $h_0$ to be $\frac{g}{\hat{\psi}} \chi_{\cup_{n \leq N} F_n}$.  Then, $h_0(\xi) \hat{\psi}(\xi) = g(\xi)$ on $\cup_{n \leq N} F_n$, so that $\|h_0 \hat{\psi} - g \| < \epsilon$.  Since $\wo$ is closed, $g \in \wo$; furthermore all such $g$'s are dense in $L^2(E)$, whence $L^2(E) \subset \wo$.

Since $W_j \perp W_k$, $2^j E \cap E$ is a set of measure zero, and since $\oplus W_j$ is dense in $\ltwo$, it follows that $\cup_j 2^j E = \Bb{R}$.  Furthermore, by theorem~\ref{T:ef}, $E$ is $2\pi$ translation congruent to $[0, 2\pi)$, hence, by theorem~\ref{T:set}, $E$ is a wavelet set, and $\psi$ is a MSF wavelet.
\end{proof}

\begin{corollary}
The equivalence class $\Cal{M}_\infty$ can be characterized in the following two ways: 
\begin{enumerate}
\item $\Cal{M}_{\infty} = \cap_{n=0}^{\infty} \Cal{L}_n$,
\item $\Cal{M}_{\infty}$ is precisely the collection of all MSF wavelets.
\end{enumerate}
\end{corollary}
\begin{proof}
By theorem~\ref{T:MSF}, $V_0$ is invariant under $\Cal{G}_n$ for all $n$ if and only if $V_0$ is invariant under translations by all real numbers.  This is equivalent to $\psi$ being a MSF wavelet.
\end{proof}

\section{A Characterization of $\Cal{M}_n$}

Suppose $\psi$ is a wavelet that is in $\Cal{L}_1$.  If $T_{\frac{m}{2}} f \in W_0$, then by taking the Fourier Transform, we have $e^{-i \frac{m}{2} \cdot} \hat{f} \in \wo$, and vice versa, so $W_0$ is invariant under translations by half integers if and only if $\wo$ is invariant under multiplication by $e^{-i \frac{m}{2} \cdot}$.  Because of this, we shall proceed with the analysis in the frequency domain.

If $f \in W_0$, then we can write $f = \sum_{k \in \Bb{Z}} c_k T^k \psi$, so taking the Fourier transform of both sides yields $\hat{f} = h \hat{\psi}$ for some $h \in L^2([0, 2\pi))$.  Hence, we can describe $\wo$ by $\{ h(\xi) \hat{\psi}(\xi) : h \in L^2([0, 2\pi)) \}$.  Suppose that $\xi \in E = supp(\hat{\psi})$.  If $\wo$ is invariant under multiplication by $e^{-i \frac{m}{2} \xi}$, then for $m = 1$,
\[ e^{-i \frac{1}{2} \xi} h(\xi) \hat{\psi}(\xi) = g(\xi) \hat{\psi}(\xi) \]
for some $g \in L^2([0, 2\pi))$.  Note that if $\xi \in supp(\hat{\psi})$, then $e^{-i \frac{1}{2} \xi} h(\xi) = g(\xi)$.  Let $\xi \in supp(\hat{\psi})$ and let $k$ be an odd integer.  Then,
\begin{align}
g(\xi) \ \hat{\psi}(\xi + 2k\pi) &= g(\xi+ 2k\pi) \ \hat{\psi}(\xi + 2k\pi) \notag \\
	&= e^{-i (\frac{1}{2}) (\xi + 2k\pi)} \ h(\xi + 2k\pi) \ \hat{\psi}(\xi + 2k\pi) \notag \\
	&= - e^{-i \frac{1}{2} \xi} \ h(\xi) \ \hat{\psi}(\xi + 2k\pi) \notag \\
	&= - g(\xi) \ \hat{\psi}(\xi + 2k\pi) \notag
\end{align}
This calculation shows that $\hat{\psi}$ cannot have both $\xi$ and $\xi + 2k\pi$ in its support.  We have established the first characterization theorem.

\begin{theorem} \label{T:L1}
Let $\psi$ be a wavelet.  Then $\psi \in \Cal{L}_1$ only if $E = supp(\hat{\psi})$ is not partially self similar with respect to any odd multiple of $2 \pi$.
\end{theorem}

\begin{corollary}
If $supp(\hat{\psi}) = \Bb{R}$, then $\psi \in \Cal{M}_0$.
\end{corollary}

\begin{corollary}
If $\psi$ is compactly supported, then $\psi \in \Cal{M}_0$.
\end{corollary}

It is interesting to note that most of the wavelets used in practice have this property.  It is unclear at this point if this has a meaningful interpretation from a numerical analysis point of view.  

Theorem~\ref{T:L1} extends to $\Cal{L}_n$ in the following natural way.
\begin{theorem} \label{T:Ln}
Let $\psi$ be a wavelet.  Then $\psi \in \Cal{L}_n$ only if the support of $\hat{\psi}$ is not partially self similar with respect to any odd multiple of $2^j \pi$ for all $j = 1,2,\ldots,n$.
\end{theorem}
\begin{proof}
Let $\psi \in \Cal{L}_n$.  Hence,
\[ e^{-i \frac{1}{2^n} \xi} h(\xi) \hat{\psi}(\xi) = g(\xi) \hat{\psi}(\xi) \]
for some $g \in L^2([0, 2\pi))$.  Let $1 \leq j \leq n$, and let $k$ be an odd integer.  Then, by a similar computation, 
\begin{align*}
g(\xi) \ \hat{\psi}(\xi + 2^j k\pi) &= g(\xi+ 2^j k\pi) \ \hat{\psi}(\xi + 2k\pi) \\
	&= e^{-i (\frac{1}{2^n}) (\xi + 2^j k\pi)} \ h(\xi + 2^j k\pi) \ \hat{\psi}(\xi + 2^j k\pi) \\
	&= e^{-i \frac{k}{2^{n-j}} \pi} e^{-i \frac{1}{2^n} \xi} \ h(\xi) \ \hat{\psi}(\xi + 2^j k\pi) \\
	&= e^{-i \frac{k}{2^{n-j}}\pi} g(\xi) \ \hat{\psi}(\xi + 2^j k\pi)
\end{align*}
as above.
\end{proof}

We have now established necessary conditions for wavelets to be in the equivalence classes $\Cal{M}_k$ for $k$ not equal to 1 or $\infty$.  This does not shed light onto whether such wavelets exist.  Fortunately, to aid in the search, the converse of theorem~\ref{T:Ln} also holds.
\begin{theorem} \label{T:Ln2}
Let $\psi$ be a wavelet and let $E = supp(\hat{\psi})$ be such that it is not partially self similar with respect to any odd multiple of $2^j \pi$ for $j=1,2,\ldots,n$.  Then $\psi \in \Cal{L}_n$.
\end{theorem}
\begin{proof}
It suffices to show that 
\[ e^{-i \frac{1}{2^n} \xi} \hat{\psi}(\xi) = g(\xi) \hat{\psi}(\xi) \] for some $g \in L^2([0, 2\pi))$.

Let $F \subset E$ be such that $\tau:F \to [0, 2\pi)$ is a bijection.  (It can be easily shown that $\tau:E \to [0, 2\pi)$ is a surjection.)  The injective property of $\tau$ can be assured in the following manner: for each $\xi \in [0, 2\pi)$, define the set $Z_{\xi} = \{m_{\xi} \in \bbz: \xi + 2 m_{\xi} \pi \in E\}$, then for $\xi$ choose $k_{\xi}$ to be 0 if $\xi \in E$, if not, choose $k_{\xi} = min\{m > 0: m in Z_{\xi}\}$, else choose $k_{\xi} = max\{m < 0: m \in Z_{\xi}\}$.  Let $F = \{\xi + 2 k_{\xi} \pi : \xi \in [0, 2\pi) \}$.  Note that by construction, $F$ is $2\pi$ translation congruent to $[0, 2\pi)$.  Hence,
\[ e^{-i \frac{1}{2^n} \xi} \chi_{F}(\xi) = g(\xi) \]
where $g(\xi) \in L^2(F)$ and is $2 \pi$ periodic.  Thus, for $\xi \in F$, 
\[ e^{-i \frac{1}{2^n} \xi} \hat{\psi}(\xi) = g(\xi) \hat{\psi}(\xi). \]

For almost any $\xi \in E \setminus F$, there exists a $\xi' \in F$ and an integer $l_{\xi}$ such that $\xi - \xi' = 2 l_{\xi} \pi$.  Moreover, by hypothesis, $l_{\xi}$ is an even multiple of $2^n$, since $E$ is not partially self similar with respect to any odd multiple of $2^j \pi$.  Since $e^{-i \frac{1}{2^n} \xi}$ is $2^n \pi$ periodic, we have that for $\xi \in E - F$,
\begin{align*}
e^{-i \frac{1}{2^n} \xi} \hat{\psi}(\xi) &=  e^{-i \frac{1}{2^n} (\xi' + 2l_{\xi}\pi)} \hat{\psi}(\xi' + 2l{\xi}\pi) \\
	&= e^{-i \frac{1}{2^n} \xi'} \hat{\psi}(\xi' + 2l{\xi}\pi) \\
	&= g(\xi') \hat{\psi}(\xi' + 2l{\xi}\pi) \\
	&= g(\xi) \hat{\psi}(\xi).
\end{align*}
This completes the proof.
\end{proof}

We have established the following characterization of the $\Cal{M}_n$'s.
\begin{corollary} \label{C:Mn}
The equivalence class $\Cal{M}_n$ consists of all wavelets $\psi$ such that the support of $\hat{\psi}$ is not partially self similar with respect to any odd multiples of $2^k \pi$, for $k=1,2,\ldots,n$ but is partially self similar with respect to some odd multiple of $2^{n+1} \pi$.
\end{corollary}

\section{Examples}

In this section, we will present examples of wavelets that are in the first four equivalence classes, with the last being in $\Cal{M}_0$ but it is not an MRA wavelet, and hence cannot be compactly supported.  The tool used to generate these wavelets is operator interpolation.  Let $\psi_{W_1}$ and $\psi_{W_2}$ be MSF wavelets, with corresponding wavelet sets $W_1$ and $W_2$, respectively.  By theorem~\ref{T:set}, $W_1$ is $2 \pi$ translation congruent to $W_2$.  If $\sigma: W_1 \to W_2$ is effected by this translation congruence, then $\sigma$ can be extended to a measurable bijection of $\Bb{R}$ by defining $\sigma(x) = 2^{-n} \sigma(2^n x)$ where $n$ is such that $2^n x \in W_1$.

If $\sigma$ is \emph{involutive}, i.e. $\sigma^2$ is the identity, and if $h_1$ and $h_2$ are measurable, essentially bounded, 2-dilation periodic functions (i.e. $h_1(2x) = h_1(x)$), then $\psi$ defined by 
\[ \hat{\psi} = h_1 \hat{\psi}_{W_1} + h_2 \hat{\psi}_{W_2} \]
is again a wavelet provided the matrix
\begin{equation} \label{E:M}
\left( \begin{matrix}
h_1 & h_2 \\ h_2 \circ \sigma^{-1} & h_1 \circ \sigma^{-1} \end{matrix}
\right) \end{equation}
is unitary almost everywhere.  (Since $\sigma^{-1}$ is 2-homogeneous, and the $h_i$'s are 2-dilation periodic, in general it suffices to check this condition on $W_1$.)  A complete discussion of this can be found in \cite{DL}.  Note that the interpolated wavelet $\psi$ has the property that $supp(\hat{\psi}) \subset W_1 \cup W_2$.  Further, note that since $\sigma$ on $W_1$ is given by translations by integral multiples of $2 \pi$, $\sigma$ completely describes the partial self similarity of $W_1 \cup W_2$ with respect to multiples of $2 \pi$.

In the following examples, $\sigma$ will always be involutive.

\begin{example} \label{E:M1}
We shall now present an example of a wavelet in $\Cal{M}_1$, which by corollary~\ref{C:Mn} is equivalent to $E = supp(\hat{\psi})$ being not partially self similar with respect to any odd multiples of $2 \pi$, but does have partially self similarity with respect to some multiple of $4 \pi$.  Consider the following two wavelet sets:

\begin{align}
W_1 =& \ [-\frac{8\pi}{7}, -\frac{4\pi}{7}) \cup [\frac{4\pi}{7}, \frac{6\pi}{7}) \cup [\frac{24\pi}{7}, \frac{32\pi}{7}) \notag \\
W_2 =& \ [-\frac{8\pi}{7}, -\frac{4\pi}{7}) \cup [\frac{2\pi}{7}, \frac{3\pi}{7}) \cup [\frac{24\pi}{7}, \frac{30\pi}{7}) \notag \\
& \ \cup  [\frac{31\pi}{7}, \frac{32\pi}{7}) \cup [\frac{60\pi}{7}, \frac{62\pi}{7}) \notag
\end{align}
A routine calculation shows:
\begin{equation}
\sigma(\xi) = \notag
  \begin{cases}
    \xi, &\text{$\xi \in W_1 \cap W_2$} \\
    \xi - 4 \pi, &\text{$\xi \in [\frac{30\pi}{7}, \frac{31\pi}{7})$} \\
    \xi + 8 \pi, &\text{$\xi \in [\frac{4\pi}{7}, \frac{6\pi}{7})$}
  \end{cases}
\end{equation}
This $\sigma$ is involutive.  Indeed, since $\sigma([\frac{30\pi}{7}, \frac{31\pi}{7})) = [\frac{2\pi}{7}, \frac{3\pi}{7})$ and $[\frac{2\pi}{7}, \frac{3\pi}{7}) = 2 [\frac{4\pi}{7}, \frac{6\pi}{7})$, for $\xi \in [\frac{30\pi}{7}, \frac{31\pi}{7})$, $\sigma^2(\xi) = \sigma(\xi - 4 \pi) = \frac{1}{2} \sigma (2(\xi - 4 \pi)) = \frac{1}{2}(2\xi - 8 \pi + 8\pi) = \xi$.  A similar computation shows that $\sigma^2$ is the identity on $[\frac{4\pi}{7}, \frac{6\pi}{7})$.

Construct $h_1$ and $h_2$ as follows:
\begin{align}
h_1 &= \chi_{W_1 \cap W_2} + \frac{1}{\sqrt{2}} \chi_{[\frac{4\pi}{7}, \frac{6\pi}{7}) \cup [\frac{30\pi}{7}, \frac{31\pi}{7})} \notag \\
h_2 &= \frac{1}{\sqrt{2}} \left( \chi_{[\frac{2\pi}{7}, \frac{3\pi}{7})} - \chi_{[\frac{60\pi}{7}, \frac{62\pi}{7})} \right) \notag
\end{align}
We need to check the condition of the matrix in equation \ref{E:M}.

It suffices to verify that the matrix is unitary on $W_1$.  Clearly, on $W_1 \cap W_2$ the matrix is unitary, indeed it is the identity there.  On $[\frac{30\pi}{7}, \frac{31\pi}{7})$, $h_1 = h_2 \circ \sigma^{-1} = \frac{1}{\sqrt{2}}$.  Furthermore, if $\xi \in [\frac{30\pi}{7}, \frac{31\pi}{7})$, $h_2(\xi) = h_2(2\xi) = - \frac{1}{\sqrt{2}}$.  Finally, $\sigma^{-1} (\xi) = \xi - 4 \pi \in [\frac{2\pi}{7}, \frac{3\pi}{7}) = \frac{1}{2}[\frac{4\pi}{7}, \frac{6\pi}{7})$, hence $h_1 \circ \sigma^{-1} (\xi) = \frac{1}{\sqrt{2}}$.  Thus, the matrix is simply:
\begin{equation*}
\left( \begin{matrix}
\frac{1}{\sqrt{2}} & \frac{1}{\sqrt{2}} \\ \frac{1}{\sqrt{2}} & - \frac{1}{\sqrt{2}} \end{matrix}
\right) \end{equation*}
which is unitary as required.  A similar computation shows that the matrix is also unitary on $[\frac{4\pi}{7}, \frac{6\pi}{7})$.
\end{example}

\begin{example} \label{E:M2}
Here we give an example of a wavelet in $\Cal{M}_2$, which by corollary~\ref{C:Mn} is equivalent to $E = supp(\hat{\psi})$ being not partially self similar with respect to any odd multiples of $2 \pi$ or $4 \pi$, but does have partially self similarity with respect to some multiple of $8 \pi$.  Consider the following two wavelet sets:

\begin{align}
W_1 =& \ [-8\pi, -\frac{112\pi}{15}) \cup [-\frac{16\pi}{15}, -\pi) \cup [-\frac{14\pi}{15}, -\frac{8\pi}{15}) \notag \\
& \ \cup [\frac{8\pi}{15}, \frac{14\pi}{15}) \cup [\pi, \frac{16\pi}{15}) \cup [\frac{112\pi}{15}, 8\pi) \notag \\
W_2 =& \ [-8\pi, -\frac{112\pi}{15}) \cup [-\frac{14\pi}{15}, -\frac{\pi}{2}) \notag \\
& \ \cup  [\frac{8\pi}{15}, \frac{14\pi}{15}) \cup [\frac{225\pi}{30}, 8\pi) \cup [\frac{224\pi}{15}, 15\pi) \notag
\end{align}
A routine calculation shows:
\begin{equation}
\sigma(\xi) = \notag
  \begin{cases}
    \xi,	 &\text{$\xi \in W_1 \cap W_2 $} \\
    \xi - 8 \pi, &\text{$\xi \in [\frac{112\pi}{15}, \frac{225\pi}{30})$} \\
    \xi + 16 \pi,&\text{$\xi \in [-\frac{16\pi}{15}, -\pi)$} \\ 
  \end{cases}
\end{equation}
As in example~\ref{E:M1}, $\sigma$ is involutive, and define $h_1$ and $h_2$ analogously:
\begin{align}
h_1 &= \chi_{W_1 \cap W_2} + \frac{1}{\sqrt{2}} \chi_{[-\frac{16\pi}{15}, -\pi) \cup [\frac{112\pi}{15}, \frac{225\pi}{30})} \notag \\
h_2 &= \frac{1}{\sqrt{2}} \left( \chi_{[-\frac{8\pi}{15}, -\frac{\pi}{2})} - \chi_{[\frac{224\pi}{15}, 15\pi)} \right) \notag
\end{align}
These functions satisfy \ref{E:M}.
\end{example}

\begin{example} \label{E:M3}
We shall now present an example of a wavelet in $\Cal{M}_3$.  Consider the following wavelet sets.

\begin{align}
W_1 =& \ [-16\pi, -\frac{480\pi}{31}) \cup [-\frac{32\pi}{31}, -\pi) \cup [-\frac{30\pi}{31}, -\frac{16\pi}{31}) \notag \\
& \ \cup [\frac{16\pi}{31}, \frac{30\pi}{31}) \cup [\pi, \frac{32\pi}{31}) \cup [\frac{480\pi}{31}, 16\pi) \notag \\
W_2 =& \ [-16\pi, -\frac{480\pi}{31}) \cup [-\frac{30\pi}{31}, -\frac{\pi}{2}) \notag \\
& \ \cup  [\frac{16\pi}{31}, \frac{30\pi}{31}) \cup [\pi, \frac{32\pi}{31}) \cup [\frac{31\pi}{2}, 16\pi) \cup [\frac{960\pi}{31}, 31\pi) \notag
\end{align}
Then, $\sigma$ is given by:
\begin{equation}
\sigma(\xi) = \notag
  \begin{cases}
    \xi, &\text{$\xi \in W_1 \cap W_2$} \\
    \xi - 16 \pi,&\text{$\xi \in [\frac{480\pi}{31}, \frac{31\pi}{2})$} \\
    \xi + 32 \pi,&\text{$\xi \in [-\frac{32\pi}{31}, -\pi)$} 
  \end{cases}
\end{equation}
Again, as in example~\ref{E:M1}, $\sigma$ is involutive; analogously define $h_1$ and $h_2$ as:
\begin{align}
h_1 &= \chi_{W_1 \cap W_2} + \frac{1}{\sqrt{2}} \chi_{[-\frac{32\pi}{31}, -\pi) \cup [\frac{480\pi}{31}, \frac{31\pi}{2})} \notag \\
h_2 &= \frac{1}{\sqrt{2}} \left( \chi_{[-\frac{16\pi}{31}, -\frac{\pi}{2})} - \chi_{[\frac{960\pi}{31}, 31\pi)} \right) \notag
\end{align}
\end{example}

\begin{example} \label{E:M0}
In this example we shall construct a non-MRA wavelet in $\Cal{M}_0$.  Consider the following wavelet sets:

\begin{align}
W_1 =& \ [-\frac{32\pi}{7}, -\frac{28\pi}{7}) \cup [-\frac{7\pi}{7}, -\frac{4\pi}{7}) \cup [\frac{4\pi}{7}, \frac{7\pi}{7}) \cup [\frac{28\pi}{7}, \frac{32\pi}{7}) \notag \\
W_2 =& \ [-\frac{8\pi}{7}, -\frac{4\pi}{7}) \cup [\frac{4\pi}{7}, \frac{6\pi}{7}) \cup [\frac{24\pi}{7}, \frac{32\pi}{7}) \notag
\end{align}
Both of these wavelets are non-MRA wavelets.  It is shown in \cite{W} that the interpolated wavelet also is not an MRA wavelet.  We have that $\sigma$ is given by: 
\begin{equation}
\sigma(\xi) = \notag
  \begin{cases}
    \xi, &\text{$\xi \in [-\frac{7\pi}{7}, -\frac{4\pi}{7}) \cup 
			[\frac{4\pi}{7}, \frac{6\pi}{7}) \cup
			[\frac{28\pi}{7}, \frac{32\pi}{7})$} \\
    \xi - 2 \pi, &\text{$\xi \in [\frac{6\pi}{7}, \frac{7\pi}{7})$} \\
    \xi + 8 \pi, &\text{$\xi \in [-\frac{32\pi}{7}, -\frac{28\pi}{7})$}
  \end{cases}
\end{equation}
Construct $h_1$ and $h_2$ as follows:
\begin{align}
h_1 &= \chi_{W_1 \cap W_2} + \frac{1}{\sqrt{2}} \chi_{[-\frac{32\pi}{7}, -\frac{28\pi}{7}) \cup [\frac{6\pi}{7}, \frac{7\pi}{7})} \notag \\
h_2 &= \frac{1}{\sqrt{2}} \left( \chi_{[-\frac{8\pi}{7}, -\frac{7\pi}{7})} - \chi_{[\frac{24\pi}{7}, \frac{28\pi}{7})} \right) \notag
\end{align}
\end{example}

% BIBLIOGRAPHY

\end{document}